\documentclass[a4paper]{amsart}

\usepackage[english]{babel}
\usepackage{amssymb}
\usepackage{MnSymbol}
\usepackage{amsmath}
\usepackage{amsthm}
\usepackage{amsaddr}
\usepackage{graphicx}
\usepackage{fancyhdr}
\usepackage{bbm}
\usepackage{multirow}
\usepackage{enumitem}
\usepackage{algorithm2e}
\usepackage{soul}
\newfont{\bssten}{cmssbx10}
\newfont{\bssnine}{cmssbx10 scaled 900}
\newfont{\bssdoz}{cmssbx10 scaled 1200}

\usepackage{latexsym,amssymb,amsthm,amsxtra}
\usepackage{amsmath}
\usepackage{stmaryrd}
\usepackage{mathrsfs}
\usepackage{tikz}
\usepackage{pgf}
\usepackage{cancel}
\usepackage{caption} 
\usepackage{dsfont}

\newtheorem{theoreme}{Theorem}

\newtheorem{definition}{Definition}

\theoremstyle{definition}
\theoremstyle{remark}

\newtheorem{proof(theorem)}{Proof(Theorem)}

\usepackage{tabularx} 
\usepackage{multirow} 
\usepackage{multicol} 
\usepackage{arydshln} 
\usepackage{fancybox} 
\usepackage{multicol} 
\usepackage{array} 
\usepackage{fancybox}
\usepackage{listings}
\usepackage{float}
\usepackage{array,multirow,makecell}
\setcellgapes{1pt}
\makegapedcells
\newcolumntype{R}[1]{>{\raggedleft\arraybackslash }b{#1}}
\newcolumntype{L}[1]{>{\raggedright\arraybackslash }b{#1}}
\newcolumntype{C}[1]{>{\centering\arraybackslash }b{#1}}
\makeatletter
\renewcommand\paragraph{\@startsection{paragraph}{4}{\z@}%
            {-2.5ex\@plus -1ex \@minus -.25ex}%
            {1.25ex \@plus .25ex}%
            {\normalfont\normalsize\bfseries}}
\makeatother
\newcommand{\E}{\mathbb{E}}

\newcommand{\R}{\mathbb{R}}

\newcommand{\x}{\mathbf{x}}
\newcommand{\y}{\mathbf{y}}
\newcommand*\diff{\mathop{}\!\mathrm{d}}

\newcommand{\Step}{\textbf{Step}}

\newcommand{\X}{\mathbb{X}}

\title[Real-time queueing models for liver transplantation]{Toward organ shortage resilient allocation policies using real-time queueing models for liver transplantation}
\author
{\small{Thomas Masanet, Benoît Audry, Christian Jacquelinet and Pascal Moyal}\\
\tiny{Université de Lorraine, Agence de la Biomédecine, Agende de la Biomédecine / Inserm U 1018 and Université de Lorraine / Inria PASTA}}

\begin{document}
\maketitle
\begin{abstract}
We report in this paper on the potential interest of real-time queueing models to optimize organ allocation policies. We especially focus on building organ shortage resilient policies in terms of equity, as we experienced differential impact of the COVID epidemic organ shortage on transplant access, according to the cause of liver failure. Patient’s death on the waiting list or dropout for being too sick, resulting from the absence of a timely available organ, is chosen as the main equity metric. 
Results obtained with the composite allocation score used in France is challenged against the so-called Early Simulated Deadline First (ESDF) real-time queueing discipline, under increasing levels of organ shortage, by extensive simulations. The ESDF policy is a variant of the well-know Earliest Deadline First (EDF) policy, which was shown as optimal in various contexts in the queueing literature. In the present case, the time to the deadline represents the remaining life duration of patients - which is of course unknown. So we propose to simulate a fictional life-duration, and give priority to the earliest simulated deadline. This leads to a simple and comprehensive representation of the system at hand by a Markov process. Our simulation results clearly show that the ESDF policy allows to maintain equity between indications, conversely to the scoring policy, which was not resilient to increasing levels of organ shortage.

\end{abstract}

\section{Introduction}
\subsection{Liver allocation in France}

Due to organ shortage, the allocation of organs for transplantation is a worldwide sensitive issue. Organs for transplantation are considered as scarce resources, as the needs exceed graft offer. For vital organs such as Heart or Liver, the absence of a timely available organ leads to the patient's adverse outcomes on the waiting list: death or dropout for being too sick for transplantation (DDTS).
Organ allocation policies try to answer to the various medical needs of transplant candidates in a suitable, efficient, equitable and traceable manner. To achieve the best use of this scarce resource, many countries moved from centre-based to patient-based allocation policies. In some cases, sequential allocation priorities related to patients conditions can be defined, as is the case for example for patients with fulminant hepatitis. In most cases, the overlaps of patients’ medical conditions do not permit to rank patients according to simple categorical criteria. In this situation, multivariate allocation scores have been demonstrated to offer efficient tools, both to support patient-based allocation policies and to optimise organ allocation in weighting the influence of various, and sometimes contradictory, equity, efficacy and feasibility allocation criteria. The Model for End-stage Liver Disease (MELD) has been demonstrated to be a good predictor of risk of death for cirrhotic patients \cite{wiesner2001meld},\cite{merion} and is widely used to prioritise liver transplant candidates. It has also been demonstrated to identify patients with an individual benefit from Liver Transplantation (LTx), and permits to avoid too early transplantation, when the risks related to the liver disease is lower than the risk of tranpslant surgical procedure \cite{merion}.

In France, the Agence de la biomédecine (ABM) is responsible for the allocation of organ and the management of the waiting lists. ABM, in collaboration with transplant professionals and patients associations defines and implement the organ allocation policy. The Liver allocation policy comprises a national ``Super-Urgency" priority for patients with fulminant hepatitis or early primary graft non function \cite{adam}. Liver are allocated to patients sharing the same blood type with the donor. In absence of SU patients,  the liver allocation score (LAS) is computed for each active patient on the national waiting list \cite{jacquelinet2008regles,calmus2008nouvelles}. LAS assigns them a priority rank and therefore allows the allocation of an organ to the patient with the highest score \cite{Handbook}. The LAS  takes into account the MELD score for cirrhotic patients (CIRRH), combined with the time spent on the waiting list for hepatocellular carcinomas (HCC) and other indications (OTHER). To handle logistic constraints,  the final score includes an interaction between the LAS and the distance between organ procurement and transplant centers using a gravity model \cite{bayer}. For some indications, referred as to ``MELD exceptions" (MXP), MELD does not reflect properly the risk of death on the waiting list. According to their condition, these patients are granted with $650$ to $800$ {additional points to the score within} $0$, $3$, $6$ \textit{or} $9$ \textit{months}, by expert advice.

\subsection{Rationale for change}
During the Covid-19 period, we observed an increase in organ shortage that resulted in an imbalance between Liver Allocation score components \cite{Annual}.
This fact, confirmed by simulation works, prompts us to find organ shortage resilient mathematical models.

\subsection{Real-time queueing systems}

{Queueing systems have become a standard model in operational research, in the performance evaluation of communication and computer networks, and supply chains, among other fields. The typical settings allow for a representation of the system at hand by a markovian stochastic process. The exact, or approximate, analysis 
of the process under study often leads to the construction of efficient tools to control, optimize or at least, simulate the system in the long run. To do so, the central question is that of the construction of efficient service disciplines (that act as a control of the system at hand) that can be implemented sequentially, 
so as to optimize a given performance parameter. In the so-called {\em real-time} context, the incoming elements have a bounded life-time, or {\em patience} in the system, and the timely execution is of critical importance. Of particular interest, in such context, is the {\em Earliest deadline first} (EDF) policy, prioritizing the customers in line having the least remaining patience time (time to the deadline) at current time. EDF policies can also model reneging, offering a mean to handle items with elapsed deadlines. 

The outcome of patients in the transplant waiting list typically stands in the field of applications of such real-time queueing systems. In the present settings, the patience time stand for the life duration of the patients, and the reneging, for their death on the waiting list, when their life duration has elapsed.}

\subsection{Outline}
In this work, we show that a practically relevant variant of the EDF policy, which we call ESDF (for {\em Earliest Simulated Deadline First}) allows to efficiently control 
the system at hand, by providing a simple representation of the system, that is amenable to extensive simulations. Second, we show that ESDF tends to be more organ shortage-resilient in terms of equity between indications, in comparison with the more classical Score policy. 

The remainder of this paper is organized as follows. In Section \ref{sec:pop}, we describe the population and study design for our simulations. Then, in Section \ref{model} we present in detail the mathematical model at hand, namely, a bipartite stochastic matching model with impatient items, as recently introduced in \cite{ADWW21,jonckheere2020generalized,masanet2022perfect,APW23}, ruled by the original ESDF matching policy. In Section \ref{sec:calibr}, we detail the calibration of the parameters of our mathematical model and in Section \ref{sec:simu}, we describe our simulation procedure. Our results are presented and analyzed in Section \ref{sec:res}. Section \ref{sec:discuss} concludes this work. Additional materials (mathematical developments and the presentation of additional algorithms in our simulation scheme) are presented in Section \ref{sec:appen}. 


\section{Population and study design}
\label{sec:pop}

\subsection{Population}

This study used aggregated data from the French national database CRISTAL. CRISTAL is a national database initiated in 1996 and administered by the ABM that prospectively collects data on all organ transplant candidates in France along with their outcomes \cite{strang}. The ABM is also in charge of the evaluation of all transplant activity including outcome of patients on the waiting list and after transplantation. Data are entered into the registry by each center and updated at least every three months. Data collection is mandatory. Withdrawal from the waiting list and patient death are prospectively notified. The study was conducted according to the French law indicating that research studies based on the national registry CRISTAL are part of the transplant assessment activity and do not require institutional reviewboard approval.
These data were used as input for the new mathematical model presented in this paper and for the survival times generation (analyses made by the ABM).


\subsection{Patients and donors arrivals/flow}
\label{flow}
\subsubsection{Patients}
All newly registered patients on the French national waiting list between January 1, 2018 and December 31, 2019, in the liver transplant centers in France were included. 3758 patients were listed for transplantation during the study period. Patients' anonymous aggregated data according to blood type, transplant indications (hepatocellular carcinoma, cirrhosis, MELD exception and others) and MELD classes to feed the mathematical model.
\subsubsection{Donors}
All brain dead deceased donors in France between January 1, 2018 and December 31, 2019, whose liver has been transplanted were included. 2458 donors were transplanted during the study period. Donors' anonymous aggregated data according to blood type were used to feed the mathematical model.

\newpage
\subsection{Survival times generation}
\label{survival}
\subsubsection{Patients}
All newly registered patients on the French national waiting list between January 1, 2012 and December 31, 2016, in the liver transplant centers in France were included. 3758 patients were listed for transplantation during the study period.
\subsubsection{Method}
A cox proportional hazard statistical was performed to predict DDTS dropout. We stratified the model by indication and only used MELD as covariate.
Inverse transformation method was used to simulate the individual patient survival time on the waiting list \cite{LSAM},\cite{bender}. 
We also used a cox model and the inverse transformation method to simulate the transition time to cirrhosis or others status to MELD exception status.

\subsection{Study design}
To mimic the impact of changes in organ shortage, we simulate the outcome of patients under various levels of liver donors. LTx and DDTS crude rates were selected as the evaluation endpoint. They were computed from the incident cohort of newly enlisted patients during the first 2 years of the study period. We'll observed the outcome of those patients on the whole 10 years study period.

We simulate 8 scenarios: 2 matching policies (described in section 3.3) and 4 levels of organ shortage defined by the decrease in the percent of liver donors available (0\%, 15\%, 30\% and 50\%).
Each scenario is performed 10 times to average the variations inherent to the random settings.

\section{Mathematical model}
\label{model}
\subsection{A matching model}
We consider a general stochastic matching model (GM), as was
defined in \cite{mairesse2016stability} (extending e.g. the classical references \cite{CKW09}, \cite{AW12} and \cite{BGM13} to single arrivals): `items' (donors and recipients) enter one by one in the system, and each of them belongs to a determinate
class. We let $V$ be the set of classes, and fix a {\em compatibility graph} $G=(V,E)$ between classes. 
Upon arrival, any incoming item of class, say, 
$i \in  V$ is either matched with an item present in the system, of a class $j$ such that $i$ shares an edge with $j$ in $G$, if any, or if no such item is available, it is stored in the system to wait for its match. 
Whenever several possible matches are possible for an incoming item $i$, a {\em matching policy} $\Phi$ determines what is the match of $i$ without ambiguity. Each matched pair (organ/recipient) departs the system right away, 
representing the transplant between the two considered item.

\subsubsection*{Recipients and donors classes} 
The set of classes is denoted by $V=D\cup R$, where $D$ is the set of donor classes, and 
$R$ is the set of recipient classes. 
For any $i\in V$, we suppose that the inter-arrival times (i.e. the durations between two consecutive arrivals) are IID (Independants and identically distributed) with the distribution $\xi_{i}$.

Each system corresponds to one single blood group, so there is a single donor class $d$, i.e.,  we set $D=\{d\}$. Each class in $R$ is characterized by a triplet $(a_1,a_2,a_3)$, where:
\begin{enumerate}
    \item $a_1$ gives the Indication of the patient. They are 4 possible indications which indicate the nature of the patient affection and which are denoted by the following terms: CIRRH for cirrhotic patients, HCC for hepatocellular carcinoma patients, MXP for patients with a MELD exception and OTHER for every other possibility. 
    \item $a_2$ is the MELD class of the patient. There are 6 different MELD classes which regroup different intervals of possible MELD value for patients:$[6,14]$, $[15,19]$, $[20,25]$, $[26,30]$, $[31,35]$ and $[36,40]$. For each indication, patients can have those 6 possible MELD Classes, except for MXP patients, who can only belong to classes $[6,14],[15,19]$ and $[20,25]$. 
    \item $a_3$ (0 or 1) indicates if the patient awaits a MELD exception. A MELD exception is given to the patients whose life quality is greatly impacted by their condition. The decision to transplant patients with MELD exception is then based on life quality, more than survival considerations. OTHER and CIRRH patients with MELD classes $[6,14],[15,19]$ and $[20,25]$ can await to receive a MELD exception. 
\end{enumerate}  

There are therefore $27$ different patient classes with $9$ OTHER classes (6 MELD classes + 3 for patients awaiting MELD exceptions), $6$ HCC classes (6 MELD classes), $9$ CIRRH classes (6 MELD classes + 3 for patients awaiting MELD exceptions), $3$ MXP classes (3 MELD classes). We denote by $R^1$ the patients who are not waiting for a MELD exception and by $R^2$ the patients awaiting a MELD exception in a way that we have $R = R_1 \sqcup R_2$.

As the patients and organs we consider all have the same blood type, all the patients in $R$ and all the donors in $D$ are compatible. Nevertheless, as patients who are awaiting a MELD exception cannot be transplanted, the set $E$ of possible matchings is given by $E = R_1 \times D$. 
In other words, the compatibility graph we consider in this work is the bipartite graph $G = (D\cup R, R_1 \times D)$.



\subsubsection*{Patience times} To stand for the timely constraint on the recipient's end, and in line with the recent extensions proposed in \cite{ADWW21}, \cite{jonckheere2020generalized}, \cite{masanet2022perfect} and \cite{APW23}, the GM model is enriched with an impatience parameter. Specifically, we suppose that a subset $R'\subseteq R$ of the classes of recipients are assigned a patience time upon arrival. If the considered item has not been matched at the end of its patience time, namely, the corresponding recipient did not find a transplant, then it leaves the system forever. The initial patience times of items of class $i\in R^{'}$ are IID, from a probability distribution $\mu_i$. Likewise, for all $i\in R'$, we assume that all elements of class $i$ are assigned predictive, or {\em simulated} patience times, which are independent of the actual patience times, and IID with the distribution 
$\mu'_i$ (which is for most classes, but not all, equal to $\mu_i$). These fictional patience times will be needed hereafter, to define a particular class of matching policies. If the predictive patience of a given recipient elapses before its match, and before reneging (i.e. before its death), then the corresponding item is still in line, and we must draw again another predictive patience time. We denote by $\gamma_{i,c}$, the distribution of predictive patience times of items of class $i$, conditional on having spent a time $c$ in the system.

\subsubsection*{Change of classes} 
Last, we suppose that our matching model with impatience is {\em dynamic}, in the sense that 
items of classes in $R$, i.e., recipients, can change classes and residual patience time within the system,  while maintaining the same place in the system. 
Specifically, we let $E_2$ be the subset of $V \times V$ gathering the set of couples $(i,j)$ such that 
a transition can be made from class $i$ to class $j$. Then, for all $(i,j) \in E_2$ we let $\lambda_{i,j}$ be the rate of transition time for a patient from class $i$ to class $j$, meaning that as soon as it becomes of class $i$ (or enter the system being of class $i$), and provided that it does not leave the queue for a transplant or a reneging before that, an item will transition to class $j$ after an exponential time (independent of everything else) of distribution $\mbox{Exp}\left(\lambda_{i,j}\right)$.
For any $i\in R'$ and $c\in\mathbb R_+^*$, we also define $\mu_{i,c}$ as the distribution of patience times of items transitioning to class $i$ after having spent a time $c$ in the system, and $\mu^{'}_{i,c}$, the generic predictive patience time of an item transitioning to class $i$ after a time $c$ spent in the system.

\subsection{Model Dynamics}
\label{dynamic1}
From the items' point of view, the typical dynamics is as follows: A recipient $r$ of class $i$ arrives in the system, say, at time $T$, and is assigned an initial waiting time set at $0$. If $i \in R^{'}$, $r$ also gets assigned a patience time $P_1$ of law $\mu_{i}$. This means that the patience time of $r$ is theoretically due to end at time $T+P_1$. The recipient $r$ also gets assigned a predictive patience $P'_1$ of law $\mu_{i}$. 

For all $j\in R$ such that $(i,j)\in E_1$, we draw an exponential random variable of law $C_{i,j}\sim \mbox{Exp}(\lambda_{i,j})$. 
Then, if one of these exponential clocks (take the smallest one, say it is $C_{i,j}$) rings before $P_1$ (if $i\in R'$), and before $r$ could find an organ for a match according to the policy $\Phi$, then, at time $T+C_{i,j}$, $r$ changes class to $j\in R$, and if $j\in R'$, $r$ is assigned a new patience time 
$P_2\sim\mu_{j,C_{i,j}}$ and a new predictive patience time 
$P'_2\sim\mu'_{j,C_{i,j}}$. In the case of patient receveing a Meld exception, the waiting time is reset at $0$, to mimic the fact that expert decisions only take into account the time since the patient has received a Meld exception. We denote by $E_3$, the subset of $E_2$ describing transitions where the waiting time is reset to $0$.

If the remaining predictive patience of an item hits $0$ after a time $c$ without any of the aforementioned events happening, and while it is of class $i$, we perform a new draw of a remaining predictive patience time, following the distribution $\gamma_{i,c}$. 

\subsubsection*{Markov chain}
The state space of the model is 
$$\X = \left(V \times \R^+ \times \overline{\R} \times \overline{\R}\right) ^*.$$
Any word $\mathbf{x} \in \X $ of length $|\mathbf x|$ can be written as 
$$\mathbf{x} =x^1x^2\,\cdots\, x^{|\mathbf x|}:= (C^1,D^1,P^1_1,P^1_2)(C^2,D^2,P^2_1,P^2_2)\,\cdots\,(C^{|\mathbf x|},D^{|\mathbf x|},P^{|\mathbf x|}_1,P^{|\mathbf x|}_2),$$ 
where for $i \in \llbracket 1,{|\mathbf x|} \rrbracket$, the variables $C^i,D^i,P^i_1$ and $P^i_2$  respectively denote the class, the waiting time in the system, the remaining patience time and the predictive patience time of the $i$-th item in the system in the order of arrivals. If the class of the $i$-th item lies outside $R$, i.e., it has no patience time and predictive patience time, then we just use the convention $P^i_1=P^i_2=+\infty$. 

The system is then updated at exponential rates when one of those two events occurs on the time interval 
$[0,h]$ (as is classically the case for exponential clocks, the probability that more than one event occurs on $[0,h]$ is a $o(h)$): 
\begin{enumerate}
    \item 
    After a time $h$, a new item $(C,0,P_1,P_2)$ enters the system, in which case we first update the state by decreasing the patience and predictive patience by $h$ and increasing the waiting time by $h$ for every item in the list. Then, any negative patience means that the corresponding item has reneged, and is deleted from the system. Therefore, we update the system to the new state 
    \begin{equation}
    \label{eq:defx'}
        \mathbf{x}' = \theta\left((C^1,D^1+h,P^1_1-h,P^1_2-h)\,\,\cdots\,\,(C^{|\mathbf x|},D^{|\mathbf x|}+h,P^{|\mathbf x|}_1-h,P^{|\mathbf x|}_2-h)\right),
    \end{equation}
    where for all $\mathbf y\in \X$, $\theta(\mathbf{y})$ is the state of $\X$ composed of the items of $\mathbf y$ having patience times in $\R_+\cup \{+\infty\}$, appearing in the same order as in 
    $\mathbf y$, and where any item $j$ whose predictive patience has reached negative values on $[0,h]$ 
    has its predictive patience re-drawn from $\gamma_{C^j,D^j+h}$. 
    Then we update to state $\mathbf{x}''$ according to the matching policy, i.e: 
    \begin{itemize}
        \item[(1a)] If the new item cannot be matched, we add the item at the lost spot in line in the system, and so $$\mathbf x''=\mathbf{x}'(C,0,P_1,P_2);$$
        \item[(1b)] Else, the new item is matched with an item in $\mathbf{x}'$, say item $x'_j$, and 
        then we just set 
        $$\mathbf{x}''=\Psi_{|j}(\mathbf{x}'),$$
        where for all $\mathbf y\in \X$, $\Psi_{|j}(\mathbf{y})$ is the word obtained from $\mathbf{y}$ by just deleting its $j$-th letter. 
    \end{itemize}
    
    \item
    After a time $h$, the $i$-th item switches from class $C^i$ to class $\tilde C^i$. 
    We first update to state $\mathbf{x}'$ defined by \eqref{eq:defx'}. 
    Then, if the $i$- th item that switches class has not reneged, we determine a new patience and a new predictive patience for this item. We perform this by drawing two random IID variables $\tilde P^i_1$ and 
    $\tilde P^i_2$ of respective laws $\mu_{\tilde C_i,D_i+h}$ and $\mu'_{\tilde C_i,D_i+h}$.
    Depending on the type of transitions we then update the waiting time of the item to 
    \[\tilde D^i=\begin{cases}
    0 &\mbox{ if $(C^i,\tilde C^i) \in E_3$;}\\
    D^i+h &\mbox{ else.}
    \end{cases}\]
    The new state of the system is then updated as $\mathbf{x}''$, obtained from by 
    substituting the letter 
    $(\tilde C^i,\tilde D^i,\tilde P^i_1,\tilde P^i_1)$ to $(C^i,D^i+h,P^i_1-h,P^i_2-h)$ in $\mathbf{x}'$. 
\end{enumerate}
This dynamics characterizes a Continuous-time Markov chain (CTMC) which we denote by $\{X_t\}$, and is reminiscent 
of the measure-valued Markov representations of EDF queues proposed in \cite{DLS01}, \cite{DM08} and \cite{Moy13}. 
For completeness, we write explicitly the infinitesimal generator of that CTMC in Section \ref{subsec:infgen}.

\subsection{Matching policies}
Similarly to \cite{mairesse2016stability,jonckheere2020generalized,masanet2022perfect}, we now formalize the notion of matching policy. 
\begin{definition}\label{def:pol}
 In a dynamic matching system with impatience, a {\em matching policy} $\Phi$ is a map from 
 $\X \times (V \times \R^+ \times \R^+ \times \R)$ to $\X$ such that for all $\mathbf x=x^1x^2\cdots x^{|\mathbf x|} \in \X,$, for all 
 $y = (C,D,P_1,P_2) \in (V \times \R^+ \times \R^+ \times \R)$, 
 \[\Phi(\mathbf x,y)=\begin{cases}
 \mathbf x y &\mbox{ if no class in $\mathbf{x}$ is compatible with }C;\\
 \Psi_{|j}(\mathbf{x}) &\mbox{ else, if the policy chooses the item corresponding to }x^j.
 \end{cases}
 \]
A matching policy $ \Phi$ is said to be {\em patience-independent}, if for all $\mathbf{x}$ and $y$,  $\Phi(\mathbf x,y)$ is independent of the patience of items in $\mathbf x$ and of the patience of $y$. 
\end{definition}

\noindent In what follows, we consider three different matching policies.

\subsubsection*{EDF}
Standing for {\em Earliest Deadline First}, the EDF policy consists of matching the incoming item with the stored compatible item that has the lowest remaining patience, if any. 

\subsubsection*{SCORE}
We say that $\Phi$ is a SCORE policy, if for all 
$$\mathbf x= (C^1,D^1,P^1_1,P^1_2)\,\cdots\,(C^{|\mathbf x|},D^{|\mathbf x|},P^{|\mathbf x|}_1,P^{|\mathbf x|}_2)\in \X,$$ the choice of $j$ in Definition \ref{def:pol} is made based on the sole knowledge of the $$S\left((C^i,D^i)\right),\,i\in \llbracket 1,|\mathbf x| \rrbracket,$$ where $S$ is a measurable map from 
$V\times \R_+$ to $\R_+$. In other words, the choice of match induced by $\Phi$ is made based on the sole knowledge of the classes and waiting times of the items in line. This policy corresponds to a simplified version of the score allocation, currently used in the France transplant system. 

\subsubsection*{ESDF}
Similarly to the EDF policy, the ESDF ({\em Earliest Simulated Deadline First}) policy consists of  matching the item in the system that has the smallest remaining predictive patience among those which are compatible with the incoming item, if any. 

\medskip

{As is well known in real-time queueing theory, the EDF policy, which prioritizes the elements that are the closest to their reneging time, is optimal in terms of feasibility \cite{Der74}, tardiness \cite{Moy08} and loss probability under various statistical assumptions (see e.g. \cite{PW88} and \cite{Moy13}). However, the very implementation of the EDF policy requires the full knowledge of the patience times of all customers in the system. In the context of organ transplants, assuming the full knowledge of patience times is clearly unrealistic, as it amounts to the knowledge of the residual lifetimes in line, for all recipients. 
Thus, we aim at using policies that do not require this full knowledge, such as SCORE or ESDF. In particular, the ESDF policy can be seen as a proxy for EDF: instead of giving priority to the recipients that have the shortest remaining lifetime (an information that is clearly not available), we emulate the lifetime by simulating a survival time upon arrival (what we call the {\em predictive} patience), from the same probability distribution. Then, priority is given to the recipient that has the smallest residual predictive patience amongst all recipients in line for which a transplant is feasible.}

\section{Calibrating the parameters}
\label{sec:calibr}

In this section, we describe the precise parameters that were used for the liver transplant model and simulation.

\subsection{Arrival rates}
The arrival rates of all classes (donors and recipients) except the Meld exceptions recipients, are calibrated by the observed arrival rates of the patients in real life, over a period of 2 years (see \ref{flow}). For Meld exception recipients, the arrival rate is null: patients first arrive into the system, awaiting a MELD exception, and only then, transition to become MXP patients (see \ref{survival}).

\subsection{Patience times}
In practice, due to organ shortage, organs are only removed from the waiting line whenever there is already a patient awaiting in line upon their arrival into the system. 
So we make the assumption that organs are immediately matched upon arrival, in other words, organs (class $\{d\}$) have no patience clock upon arrival. 

The patience times of patients of all classes are calibrated by the observed death of patients of such classes in the waiting line, during a period of 4 years (see \ref{survival}). Transplanted recipients are  censored on their transplant date, since we could not observe their death while they were in the waiting list {(this is the so-called {\em censored at transplant})}. The patience time distribution is calibrated as follows.  
\begin{definition}
    We assume that the patience times of classes in $R'$ are IID from a Cox model. In other words, we have the following expression for the {\em instant risk of death} at time $t$, formally defined as the probability of death in an infinitesimal amount of time after $t$, conditional on the patient being still alive at time $t$:
$$\lambda(t,\mbox{{\em \tiny{MELD}} }) = \lambda_0(t) \exp(\sum_{i=1}^6\beta_i X_i),$$
where $\lambda_0(t)$ is the baseline instant risk of death at $t$, 
$$X_i =\mathds{1}_{\lbrace \mbox{{\em \tiny{MELD}} } \in  \mbox{ {\em \tiny{category $i$}}} \rbrace},\, i\in \llbracket 1 , 6 \rrbracket $$ 
are the categorial covariable corresponding to the MELD categories (among $[6,14]$, $[15,19]$, $[20,25]$, $[26,30]$, $[31,35]$ and $[36,40]$), and for all $i$, $\beta_i$ is a weight parameter for the categorial covariable $X_i$.
\end{definition}

\noindent Then, there are three cases for patients survival laws:
\begin{enumerate}
\item[(i)]
The real and predicted patience times of patients who are not MXP and do not await a MELD exception have the same law, based on that Cox model.
\item [(ii)]
MXP patients have a different law for their real and predicted patience. Indeed, for the MXP patients the predicted survival time that we may want to use to decide for a transplantation does not make much sense, since we mainly transplant them based on their life quality. Therefore we calibrated the predictive patience of MXP patients on the transplants observed in the waiting queue over a period of 4 years, see \ref{survival}. This guarantees that transplants for MXP patients that were based on predictive patience would not be too far from what we would expect in reality.
    
\item [(iii)]
Patients who awaits a MELD exception are the only patients who do not have a real or predicted patience time (these correspond to the only classes in $R\setminus R'$). Indeed, in our model we consider that we know from the beginning that these patients will receive a MELD exception, and have to wait for it (while in real life, the decision to give a MELD exception and the receiving of such a component is simultaneous).  
This difference makes it so that if patients in our model could renege while awaiting a MELD exception, we might not have the right rates of arrival for MXP patients. Therefore, we consider that those patients cannot renege or be transplanted while they await for a MELD exception. Clearly, this is just equivalent to saying that those patients only arrive in the system when they receive their Meld exception, {because in the context of organ shortage, under both SCORE and ESDF they would not have access to a transplant before receiving their Meld exception.} 
\end{enumerate}

The patience times described above are assigned to patients upon their arrival in the system. We will describe in Section \ref{transition} how we calibrate these patience times after a class transition. Before that, we need to describe how we re-draw the predictive patience times when the latter become negative:

\begin{definition}
    Let $i \in V^{'}$ and $c>0$. Then the law for re-drawing the predictive patience time of a recipient whose predictive patience has elapsed after having spent $c$ units of time in the system, is given by 
    \begin{equation}
    \label{eq:defgammaic}
        \gamma_{i,c}(p) = \mu{'}_{i\vert \geq c}(p) - c,\quad p\ge c,
    \end{equation} 
    where $\gamma_{i,c} = \mu'_{i\vert \geq c}$ is the following probability distribution: for all $p\ge c$, 
\begin{equation*}
\mu{'}_{i,\vert \geq c}(p)={\mu'_i(p)\over \mu'_i([c,\infty))}=\mathbb P(P_2=p \mid P_2 \ge c),
\end{equation*}
for $P_2$ a generic random variable of law $\mu'_i$. 
\end{definition}

\subsection{Class transitions}
\label{transition}
We model 2 kinds of class transitions:
\begin{enumerate}
    \item
Transition between `patient awaiting a MELD exception' to `MXP patient';
    \item
Transitions between MELD classes for patients who are not MXP nor awaiting a MELD exception, and don't have a contra-indication.
\end{enumerate}


\medskip

The transitions to receive a MELD exception are based on the observed transition over a period of $2$ years, and also take the form of a Cox Model. These transitions are the only transitions in $E_3$, i.e., we reset the time spent in the system when we give a MELD exception. For those transitions to a class $i$ after having spent $c$ units of time in the system, we draw the new real patience from the distribution 
\[\mu_{i,c} = \mu_{i\vert \ge c} -c, \] 
associated to $\mu_i$ as in \eqref{eq:defgammaic}. Then, the predictive patience is drawn from $\mu^{'}_i$, which makes sense because the predictive patience distribution was established using the observed transplantation only after a patient received a MELD exception.

The transitions between MELD classes are the transitions for which we had the least amount of data, so we calibrated the rates in the following way:
For any $i\in R$, we denote by $W^{up}_i$ (resp. $W^{down}_i$) the set of classes having the same indication, contra-indication status and MELD exception awaiting status than $i$, and having a higher (resp. lower) MELD.
\begin{enumerate}
    \item[-]
    If $W^{down}_i$ is empty, then for all $j \in W^{up}_i$ we set 
    $$\lambda_{i,j} = \frac{1}{2}\frac{\xi_{j}}{\dsum_{k \in W^{up}_i} \xi_{k}}\cdot$$
    \item[-]
        If $W^{up}_i$ is empty, then for all $j \in W^{down}_i$ we let $$\lambda_{i,j} = \frac{1}{2}\frac{\xi_{j}}{\dsum_{k \in W^{up}_i} \xi_{k}}\cdot$$
    \item[-]
    Otherwise, for all $j \in W^{up}_i$ we set
    $$\lambda_{i,j} = \frac{1}{2}\frac{1}{3}\frac{\xi_{j}}{\dsum_{k \in W^{up}_i} \xi_{k}}=\frac{1}{6}\frac{\xi_{j}}{\dsum_{k \in W^{up}_i} \xi_{k}},$$
    and for all $j \in W^{down}_i$, 
        $$\lambda_{i,j} = \frac{1}{2}\frac{2}{3}\frac{\xi_{j}}{\dsum_{k \in W^{up}_i} \xi_{k}}=\frac{1}{3}\frac{\xi_{j}}{\dsum_{k \in W^{up}_i} \xi_{k}}\cdot$$
\end{enumerate}

To better understand this calibration choice, note that setting these rates for MELD transitions is equivalent to applying the following procedure:
\begin{itemize}
\item[1)]
Draw an exponential of parameter $\lambda_i = \frac{1}{2}$ that determines the time until the patient changes MELD class (whatever the MELD class is), 
setting at $2$ years, the average time for changing MELD class.
\item[2)]
If the patient MELD can decrease (i.e., $W^{down}_i$ is not empty) and/or decrease (i.e., $W^{up}_i$ is not empty) then with probability $\frac{1}{3}$ the patient gets a lower MELD and with probability $\frac{2}{3}$, a higher MELD. This conveys the idea that a patient status is more likely to deteriorate than to 
improve.
\item[3)] 
Once we have chosen if the patient status improves or deteriorates (say it deteriorates), we pick the class $j$ with probability $\frac{\xi_{j}}{\dsum_{k \in W^{up}_i} \xi_{k}}$. This means that we consider picking a new MELD class among the set $W^{up}_i$ of possible ones, which is equivalent to seeing 
the arrival of a new patient, conditional on her class belonging to $W^{up}_i$.

For those transitions to a class $j$ after a time $c$ spent in the system, we draw the new actual and predicted patience times respectively from the (conditional, and shifted) laws $\mu_{j,c}$ and $\mu^{'}_{j,c}$, respectively associated to $\mu_j$ and $\mu'_j$ respectively, as in \eqref{eq:defgammaic}. 
\end{itemize}

\subsection{Adapting the model to the discrete-time settings}
The model introduced so far is formally defined in continuous time. To adapt the present model to discrete simulations, we made the following changes and approximations:
\begin{itemize}
    \item
  Arrivals occur at each time step, and the class of the arriving patient/organ is drawn according to the initial inter-arrival rates $\xi_i$, 
  by uniformization.

    \item The law of transitions between MELD classes is geometric, with the same mean as the exponential laws of the MELD transitions.
    \item 
    For practical reasons, instead of considering that a patient awaiting a MELD exception has no patience time, we draw a predicted and real patience 
    times conditioned to be greater than the transition time to receive a MELD exception, thereby ensuring that the corresponding patients will not die before receiving the MELD exception.
\end{itemize}

\section{Simulation procedure}
\label{sec:simu}

By conducting this simulations, our aim is to observe how different policy decisions (SCORE, ESDF) and different organ shortage percentages 
(0\%, 15\%, 30\%, 50\%) affect the evolution of a patient/organ waiting queue. To do so, we split the simulation into 2 phases:

\subsection{The initiation phase}
In this phase, our aim is to create a queue that will be used as an initial state for the study of the impacts of policy decisions and organ shortage in the study phase. To do so, we start with an empty queue, and then simulate the dynamics of the queue during a reasonable amount of time 
(here, approximately 15 years), by {\em not} applying the different organ shortage percentages during this phase. 
Indeed, if it was not so and if we started to study patients and organs from an empty queue, then it would lead to an unrealistic view of the queue dynamics, as most of the first patients would be transplanted, and there would not be enough competition between patients. 

Notice that, in a standard Markovian model, one typically studies the system behavior from any possible starting position, so that it can be asked why the initial queue generated in this initiation phase is not simply viewed as an additional parameter for the study phase. 
The fact is that the composition of a patient/organ waiting line is an intricate of the positions of patients/organs in the queue, of the time they already spent in the queue, their class, etc., so that viewing the starting queue as an additional parameter would require too many additional simulations. 

\subsection{The study phase}

Starting from the queue generated in the initiation phase, we then simulate the dynamics of the queue, by applying the fixed matching policy 
and the given organ shortage percentage, during a period of approximately $10$ years. 
While performing these simulations, we record different characteristics composing the fate of organs and patients during this period:

\begin{itemize}
    \item
 The initial class they had when entering the queue (if they received a MELD exception changing their class to a new class, then we consider this new class as their initial class);
    \item
If they were transplanted, deceased or are still alive at the end of the study phase;  
    \item
The time they have spent in the system during the study phase.     
\end{itemize}

It is standard in the medical context to only study the fate of patients and organs that arrived during the study phase. However, 
to have a more complete view, we also recorded the fate of patients/organs from the initiation phase. 

\subsection{Simulations dynamics}
The simulations principle is the same whether we are in the initiation phase or the study phase: 
We divide the phase in equal amount of times called {\em time steps}. At each time step, we randomly draw an incoming patient/organ class according to a probability law $P_{arr}$ (as mentioned in \ref{flow}, $P_{arr}$ is based on the observed incoming rates of patients/organs of different classes during a period of $2$ years).

If it is indeed a patient, we also draw 
if she awaits a MELD exception. We then draw her real survival time using the aforementioned Cox model, as well as the predicted survival time in the case of the ESDF policy. If it is an organ, we simply consider a constant survival time since, from a certain time point, all organs are immediately transplanted.
We then actualize the status of every patient/organ in line, as follows: 
\begin{enumerate}
    \item [1)] Their elapsed sojourn time in the system increases by one time unit;
    \item[2)] Their real and predicted survival time, as well as the eventual time they have to wait to receive a MELD exception, decrease by one; 
    \item[3)] We determine if the patient/organ is deceased, in which case it is removed from the queue; 
    \item[4)] We determine if the patient/organ receives a MELD exception, in which case we change its class, set its elapsed sojourn time to $0$, and 
    then re-draw a real and predicted survival time conditioned 
    to be larger than the time she already spent in the system.
    \item[5)] We determine if the patient/organ state deteriorates or improves, in which case we change its MELD class.
\end{enumerate} 
We then check, whether the incoming patient/organ can be matched with an already existing organ/patient in line. If so, we determine its match 
according to the matching policy, and we then remove the chosen organ/patient from the system. If not, we add the incoming patient/organ to the queue.

We then repeat this operation for every time step. 
Note that we can separate the phase in steps of time because the inter-arrival times of our model are supposed constant. 
At the end of the two phases, we return the state of the queue for the initiation phase, and the fate of all arrived elements during the study phase. 

\subsection{Main algorithms}

\begin{figure}[h!]
\begin{center}
\includegraphics[scale=0.55]{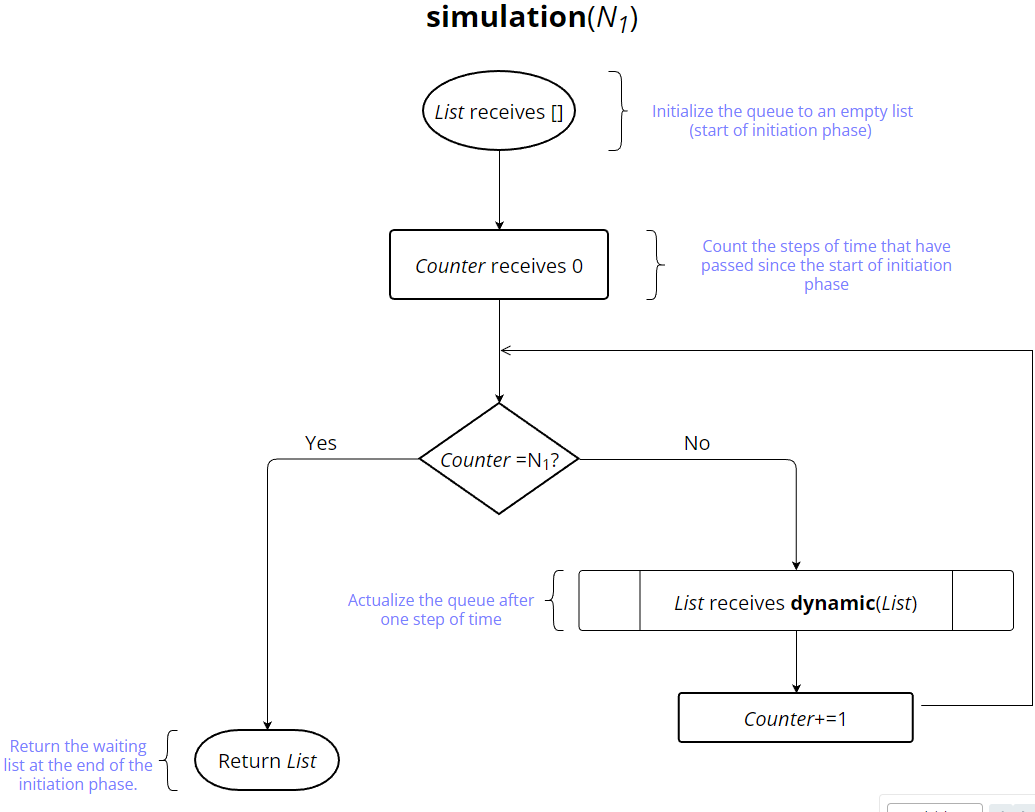}
\caption[smallcaption]{\textbf{simulation} function }
\label{fig:simulation}
\end{center}
\end{figure}

The \textbf{simulation} function that is represented in Figure \ref{fig:simulation} is the main algorithm, with $N_1$ being the number of steps we want to simulate the queue in the initiation phase starting from an empty queue.

\begin{figure}[h!]
\begin{center}
\includegraphics[scale=0.55]{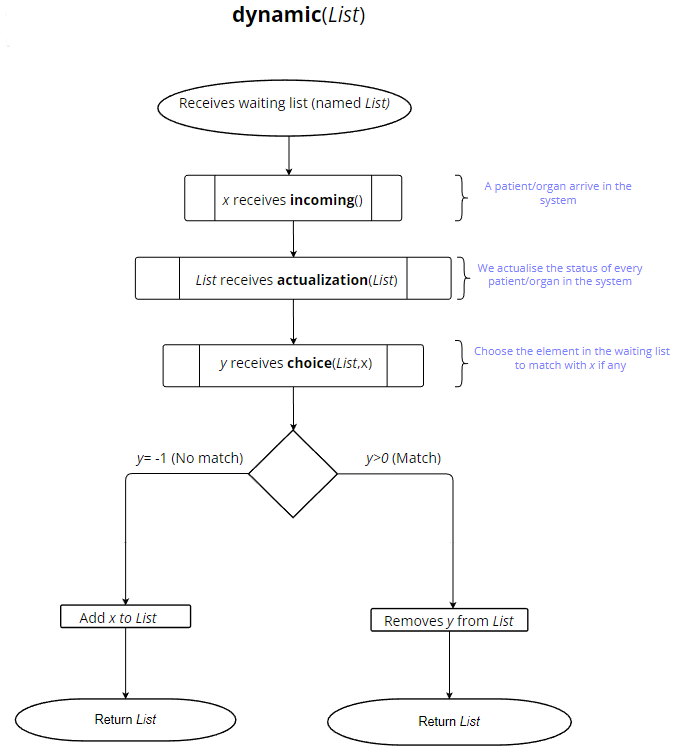}
\caption[smallcaption]{  \textbf{dynamic} function}
\label{fig:dynamic}
\end{center}
\end{figure}

The \textbf{dynamic} function represented in Figure \ref{fig:dynamic} describes the dynamics of one simulation step.

The secondary algorithms that are used in the simulation are described in section \ref{complementary}.

\section{Results}
\label{sec:res}
\subsection{Patient population}

\begin{table}[h!]
\label{Table_Incide}
\begin{center}
\includegraphics[scale=0.55]{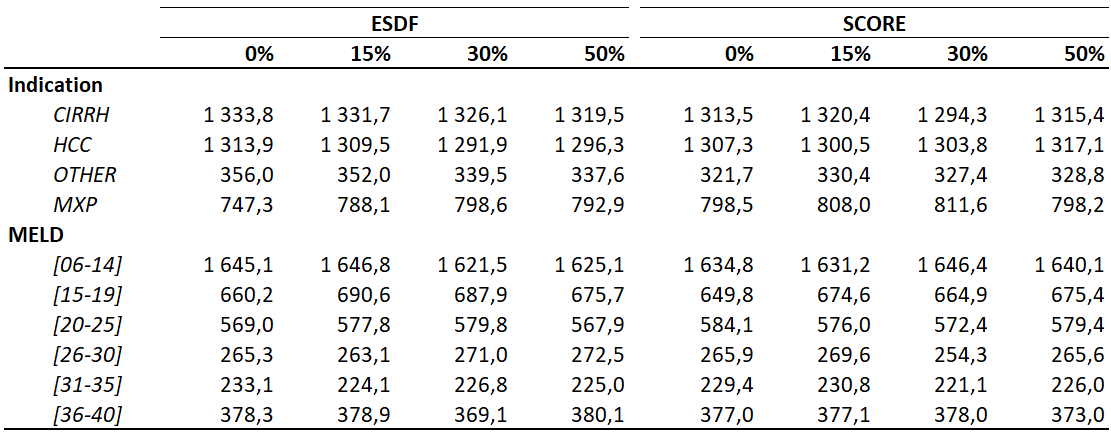}
\caption[smallcaption]{Description of the patient population (mean number of patients per 10 simulations)}
\label{Table_Incid}
\end{center}
\end{table}

The mean numbers of patients enlisted during the first 2 years of the study period per 10 simulations are presented in Table \ref{Table_Incid}. On columns we have the 2 matching policies (ESDF and SCORE) and the 4 levels of organ shortage (0\%, 15\%, 30\% and 50\%) which correspond to all the simulated scenarios. On rows we have the 4 indications and the 6 MELD classes. 

For instance, the first column for ESDF policy and 0\% of organ shortage show that the mean number of patient arrivals for CIRRH indication is 1333.8 and for [36-40] MELD class is 378.3. 

Variations in patient arrival flow by matching policies and organ shortage scenarios due to the stochastic modeling are low, and are likely to have a very small influence on DDTS and LTx rates between the two policies.

CIRRH and HCC are the 2 main indications with an average of nearly 1300 patients in 2 years. OTHER is the smallest one with an average of 321 to 356 patients and MXP has an average of 747 to 811 patients. As expected, according to the actual epidemiology of liver diseases, the  number of patients is decreasing with increasing MELD class.

\subsection{Ltx and DDTS crude rates}


\begin{figure}[h!]
\begin{center}
\includegraphics[scale=0.55]{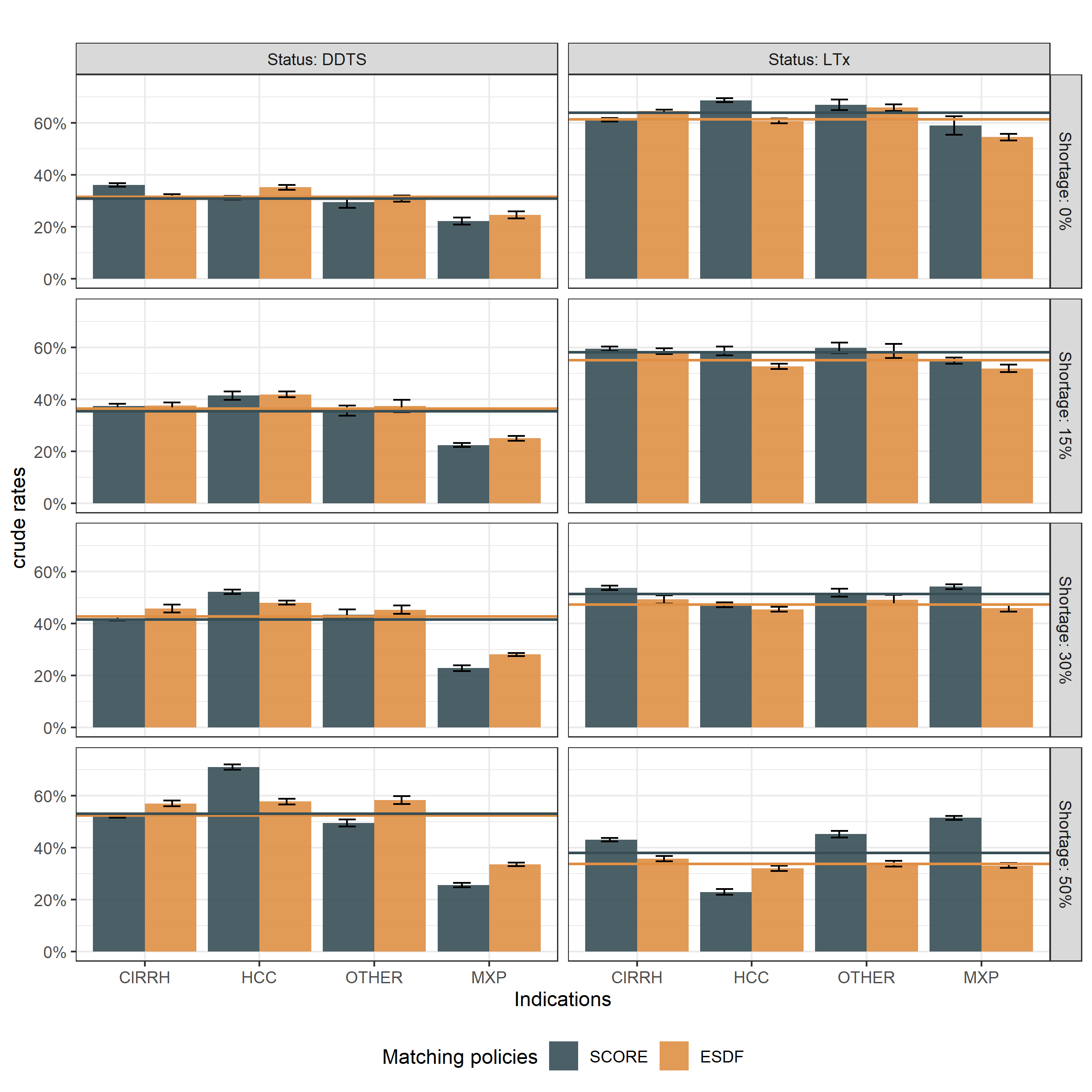}
\caption[smallcaption]{ DDTS and LTx rates by organ shortage level and matching policy. }
\label{DDTS}
\end{center}
\end{figure}

Figure \ref{DDTS} reports on the crude DDTS and LTx rates by organ shortage level and matching policy. SCORE and ESDF are displayed respectively in darkgray and orange bars. 
Of note, as we stated that an equitable allocation system should share the burden of organ shortage independently of patients characteristics, we defined DDTS variance as a relevant judgement criterion for indications where DDTS is at stake (CIRRH, HCC, OTHER).
MXP patients whose LTx indication is mainly driven by a poor quality of life have of course the most low DDTS rates. For this indication, the objective is to provide fair LTx rates.

As expected, increasing organ shortage is associated with higher DDTS rates (32\% to 52\% for ESDF policy and 31\% to 53\% for SCORE policy) and lower LTx rates (61\% to 34\% for ESDF policy and 64\% to 38\% for SCORE policy). 

ESDF and SCORE policies have roughly equivalent  DDTS rates (no significant differences) regardless of the level of shortage. But the main point of this figure is that the DDTS rates for CIRRH, HCC and OTHER indications remain equitable with increasing organ shortage level With ESDF policy, in contrast to SCORE policy where DDTS rates for HCC increase significantly in case of 30\% and 50\% increase in organ shortage. 

The same holds true for LTx rates which remain equitable for all indications with the ESDF polcy.

In contrast, the SCORE policy provides inequitable LTx rates with higher LTx rates for MXP as compared to CIRRH, OTHER and HCC, the later having the lowest LTX rate that might explain the inequitable increase of DDTS for this indication.

\subsection{Variance DDTS rate}


\begin{figure}[h!]
\begin{center}
\includegraphics[scale=0.55]{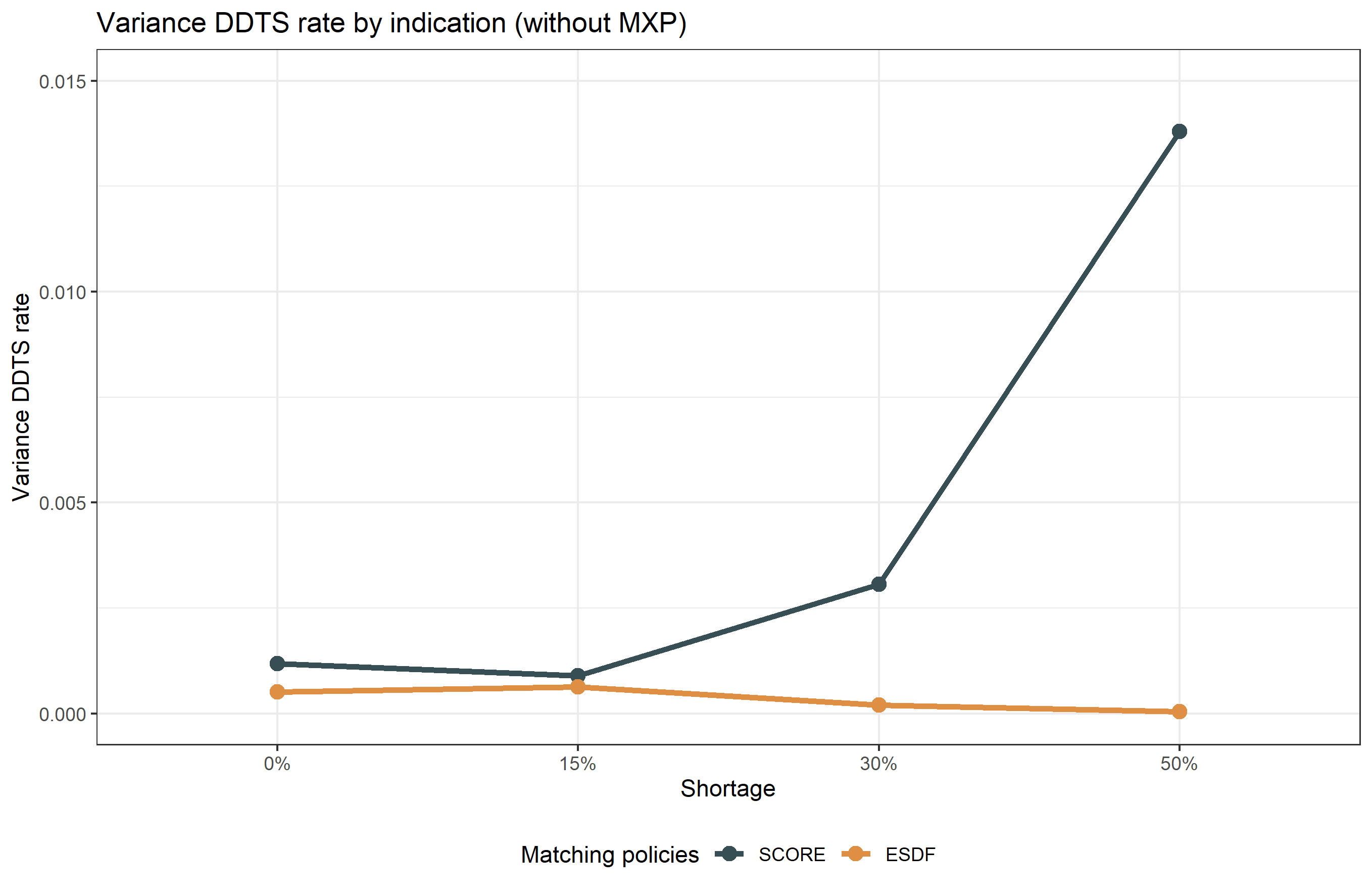}
\caption[smallcaption]{Variance DDTS rate by indication (without MXP) }
\label{variance}
\end{center}
\end{figure}
As Figure \ref{variance} clearly demonstrates, under the ESDF policy, the variance of DDTS rates for CIRRH, HCC and OTHER indications remains close to 0 for all organ shortage levels.
In contrast, the variance of DDTS increases with the organ shortage for the SCORE policy.

This result demonstrates that the ESDF policy outperforms the SCORE policy on our main equity criterion based on the DDTS rates for CIRHH, HCC and OTHER, while providing equitable LTx rates for all indications.

\section{Discussion}
\label{sec:discuss}

\subsection{Main findings}
Of course, DDTS rates will always increase and LTX rates always decrease with organ shortage, but a major result of our study is to demonstrate that ESDF policy offers prospects for allocation system resilient to organ shortage in terms of equity on both DDTS and LTx rates. 
Of course, increasing organ shortage in liver transplantation leads to an increase of DDTS rates. 
Of note, at each level of organ shortage, ESDF policy performs similarly to the Score policy in terms of DDTS rate. 

Moving from a MELD-based SCORE to a continuous expected lifetime ESDF policy (based on the projected remaining patience into the system) seems like a game changer. This study proposes an extension to the classical model EDF to a more realistic {\em Earliest Simulated Deadline First} (ESDF) policy, allowing to emulate the behaviour of the EDF policy, by predicting the value of the remaining patience times of the patients. 

In our study, organ shortage levels were related to decrease in organ procurement. Nevertheless, our results are likely valid for epidemiological or clinical situations resulting in an increase in LTx indications.

\subsection{Limitations and Future works}

Despite the good results obtained in this study, there are some limitations.

Of note, (i) the initiation phase provided a higher number of transplant candidates for the study phase, resulting in a higher organ shortage in this study than the actual situation in France ; (ii) we used 10 years DDTS rates whereas 2 or 5 years DDTS rates are usually considered.  Anyway, this does not affect the key results related to the resilience of ESDF policy to organ shortage.
Last, while we hope that a 30\% decrease in organ procurement won't be observed in practice, this corresponds to a Organ/LTx candidate ratio that might occur in some countries. Outside such a decrease, a MELD based SCORE system might be resilient to mild decreases in organ procurement.

This study uses a macro-population simulation Markov model with classes.
Future works is to validate these promising results using our micro-population simulation platform.
\subsection{Conclusion/Take home message}
Continuous expected lifetime ESDF policy might provide organ shortage resilient allocation systems in terms of equity, without significant impact on the overall DDTS rate.

\newpage

\section{Supplementary materials}
\label{sec:appen}

\subsection{Infinitesimal Generator}
\label{subsec:infgen}
The dynamics of the Markov chain $\{X_t\}$ as presented in \ref{dynamic1} can also be  encoded by the following Infinitesimal generator. 
Hereafter, fix 
$$\mathbf y = y_1, \cdots, y_{\vert \y \vert} =(C^1,D^1,P^1_1,P^1_2)(C^2,D^2,P^2_1,P^2_2)\,\cdots\,(C^{|\mathbf y|},D^{|\mathbf y|},P^{|\mathbf y|}_1,P^{|\mathbf y|}_2) \in \X,$$
and let us introduce the following additional piece of notation, 
\begin{enumerate}
\item \underline{Aging:}\\\\
To represent the aging of items, for all $h>0$ we set 
\begin{multline*}
\tau_h(y) = (C^1,D^1+h,P^1_1-h,P^1_2-h)\,\cdots\,(C^{|\mathbf y|},D^{|\mathbf y|}+h,P^{|\mathbf y|}_1-h,P^{|\mathbf y|}_2-h).
\end{multline*}
\medskip
\item \underline{Removing reneged recipients and redrawing predictive patience times:}\\\\
    Let $O^{\y}\lbrace i_1,\cdots, i_l\rbrace $, with $i_1 <i_2<\cdots < i_{l}$, be the set of indices of those items of $\y$ having positive patience and negative predictive patience.
    Then, for all $ (p{'}_1,\cdots,p^{'}_{l}) \in \R_+^l$, we let 
    $\theta_2(\y,(p{'}_1,\cdots,p^{'}_{l}))$ be the state of $\X$ composed of the items of $\mathbf y$ having patience times in $\R_+\cup \{+\infty\}$, appearing in the same order as in $\mathbf y$ and where the item in position $i_k$ has a new predictive patience $p_k$, for all $k \in \llbracket 1, l\rrbracket$.
\bigskip
\item \underline{Choosing and removing matched items:}\\\\
Let $U^{\y}$ be the subset of $V$ of those classes that can be matched with elements represented by $\y$. Then, 
for any admissible matching policy $\Phi$ and any class $i$, we let $I(\Phi,\mathbf{y},i)$ be the index of the item of $\mathbf{y}$ that will be matched with an incoming item of class $i$, according to $\Phi$.


\bigskip
\item \underline{Transitioning items:}\\\\
Let $E_2(\y)$ be the set of possible transitions of items of $\y$, where any transition for an item $I$ is denoted $(i,j,k)$ where $i$ is the class of $I$, $j$ is a class available for $I$ to transition, and $k$ is the position of $I$ in $y$. Similarly, we let $E_3(\y)$ be the set of possible transitions of items of $\y$ where the timer is reseted. We then define:
\begin{align*}
 E_4(\y) &= \lbrace (i,j,k) \in E_2(\y) \backslash E_3(\y), j\in R^{'} \rbrace,\\
 E_5(\y) &= \lbrace (i,j,k) \in E_3(\y), j\in R^{'} \rbrace,\\
 E_6(\y) &= \lbrace (i,j,k) \in E_2(\y) \backslash E_3(\y), j\in R\backslash R^{'} \rbrace,\\
 E_7(\y) &= \lbrace (i,j,k) \in E_3(\y) , j\in R\backslash R^{'} \rbrace.
 \end{align*}
Let $k \in \llbracket 1, \vert y \vert \rrbracket$ and $(C,D,P,P^{'}) \in V\times R^+\times \R \times \R$. Then $\theta_3(\y,k,(C,D,P,P^{'}))$ is the state of $\X$ composed of the items of $\mathbf y$, where the item in $k$-th position is replaced by $(C,D,P,P^{'})$.
\end{enumerate}
Gathering the arguments of section \ref{dynamic1}, we obtain that 
\begin{theoreme}
The infinitesimal generator $\mathscr A$ of the Markov chain $X$ is characterized as follows: 
for all sufficiently smooth functions $F: \X \mapsto \R$, for all $\x \in \X$ such that $O^{\x} = \lbrace i_1,\cdots, i_{l} \rbrace $  with $i_1 <i_2<\cdots < i_{l}$, we have that 
\scriptsize{\begin{multline*}
    \mathscr AF(\mathbf{x})  = \lim_{h \rightarrow 0} \frac{\E[F(X_{k+h})-F(X_k) \vert X_k = \mathbf{x}]}{h}\\
    = \lim_{h \rightarrow 0} \sum_{\ell\in U_{\mathbf{x}}} \xi_\ell \int_{\R_+^{l}} \bigl[F(\Psi_{|I(\Phi,\theta_2(\tau_h(\x),p^{'}_1,\cdots,p^{'}_{l}),i)}(\theta_2(\tau_h(\x),p^{'}_1,\cdots,p^{'}_{l}))-F(\x) \bigr] \otimes_{m=1}^l\diff \gamma_{C_{i_{m}},D_{i_m}}(p^{'}_m)\\
    + 
    \sum_{\ell\in R^{'} \backslash U_{\mathbf{x}}} \xi_\ell \int_{\R_+^{l}} \int_{\R_+^2} \bigl[F(\theta_2(\tau_h(\x),p^{'}_1,\cdots,p^{'}_{l})(i,0,p,p^{'}))-F(\x) \bigr] 
    \otimes_{m=1}^l\diff \gamma_{C_{i_{m}},D_{i_m}}(p^{'}_m) \diff \mu_{i}(p) \diff \mu_{i}^{'}(p^{'}) \\
    \shoveleft{+  \sum_{\ell\in (V\backslash R^{'})\backslash U_{\mathbf{x}}} \xi_\ell \int_{\R_+^{l}}\bigl[F(\theta_2(\tau_h(\x),p^{'}_1,\cdots,p^{'}_{l})(i,0))-F(\x) \bigr] 
    \otimes_{m=1}^l\diff \gamma_{C_{i_{m}},D_{i_m}}(p^{'}_m)}\\
    \shoveleft{+ \sum_{(i,j,k) \in E_4(\x)} \lambda_{i,j} \sum_{\ell\in R^{'} \backslash U_{\mathbf{x}}} \xi_\ell \int_{\R_+^{l}} \int_{\R_+^2} \bigl[F(\theta_3(\theta_2(\tau_h(\x),p^{'}_1,\cdots,p^{'}_{l}),k,(j,D_k,p,p^{'}))-F(\x) \bigr]}\\ \shoveright{\otimes_{m=1}^l\diff \gamma_{C_{i_{m}},D_{i_m}}(p^{'}_m)\diff \mu_{j,D_k}(p) \diff \mu_{j,D_k}^{'}(p^{'})}\\
    \shoveleft{+  \sum_{(i,j,k) \in E_5(\x)} \lambda_{i,j} \sum_{\ell\in R^{'} \backslash U_{\mathbf{x}}} \xi_\ell \int_{\R_+^{l}} \int_{\R_+^2} \bigl[F(\theta_3(\theta_2(\tau_h(\x),p^{'}_1,\cdots,p^{'}_{l}),k,(j,D_k,p,p^{'}))-F(\x) \bigr]}\\ \shoveright{\otimes_{m=1}^l\diff \gamma_{C_{i_{m}},D_{i_m}}(p^{'}_m)\diff \mu_{j,D_k}(p) \diff \mu_{j,D_k}^{'}(p^{'})}\\
    +  \sum_{(i,j,k) \in E_6(\x)} \lambda_{i,j} \sum_{\ell\in R^{'} \backslash U_{\mathbf{x}}} \xi_\ell \int_{\R_+^{l}}\bigl[F(\theta_3(\theta_2(\tau_h(\x),p^{'}_1,\cdots,p^{'}_{l}),k,(j,D_k))-F(\x) \bigr] \otimes_{m=1}^l\diff \gamma_{C_{i_{m}},D_{i_m}}(p^{'}_m) \\
    +  \sum_{(i,j,k) \in E_7(\x) } \lambda_{i,j} \sum_{\ell\in R^{'} \backslash U_{\mathbf{x}}} \xi_\ell \int_{\R_+^{l}} \bigl[F(\theta_3(\theta_2(\tau_h(\x),p^{'}_1,\cdots,p^{'}_{l}),k,(j,D_k))-F(\x) \bigr] 
    \otimes_{m=1}^l\diff \gamma_{C_{i_{m}},D_{i_m}}(p^{'}_m) \\
    + \left\{\frac{1-h\left(\sum_{i\in V} \xi_i- \sum_{(i,j,k) \in E_2(\x)}\lambda_{i,j}\right)}{h}\right\} \int_{\R_+^{l}} \bigl[F(\theta_2(\tau_h(\x),p^{'}_1,\cdots,p^{'}_{l}))-F(\x) \bigr] 
    \otimes_{m=1}^l\diff \gamma_{C_{i_{m}},D_{i_m}}(p^{'}_m).
\end{multline*}}
\end{theoreme}
\subsection{Secondary algorithms}
\label{complementary}
\subsubsection{Notations}
Every element (patient or organ) $i$ in the system can be tracked by a sextuple $[a,b_0,b_1,c,d,e]$ called the {\em information} on $i$, where:
\begin{enumerate}
    \item 
    $a$ is the class of $i$. It belongs to $R$, as defined in section \ref{model}; 
    \item 
$b_0 $ and $b_1$ are the real and predicted remaining survival time of $i$; 
    \item 
$c \in R^+$ is the time already spent in line by the patient/organ. 
    \item $d$ determines if the patient awaits to receive a MELD exception, and if so, in how much time. If $e\leq - 1$ then the patient doesn't await a MELD exception. If $-1 < d \leq 0$ the patient must receive a MELD exception. (As we decrease $e$ in a discrete way, we need a margin away from $0$ 
    that is higher than the time step.) If $d>0$, then $d$ is the time until the patient receives a MELD exception.
    \item 
$e$ only intervenes in the study phase and is a number that identify the patient/organ, negative for prevalent patients, positive otherwise.
\end{enumerate}

\subsubsection{Algorithms}

\paragraph*{Organ/patient arrival}

\begin{figure}[h!]
\begin{center}
\includegraphics[scale=0.55]{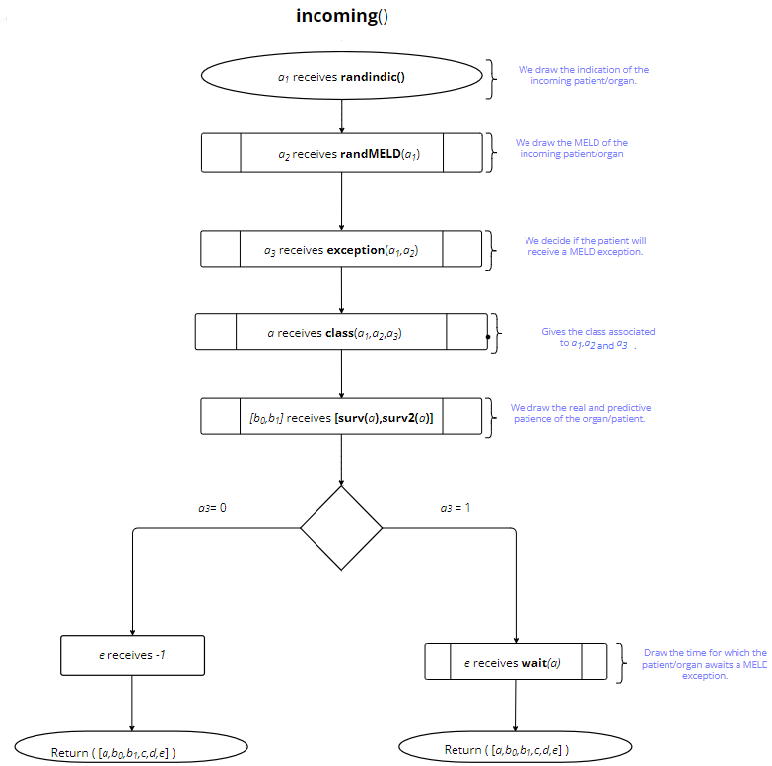}
\caption[smallcaption]{\textbf{incoming} function }
\label{incoming}
\end{center}
\end{figure}

The program shown in Figure \ref{incoming} describes the arrival of an organ/patient at one time step. It draws the information of the organ/patients, i.e their class, their eventual temporary contra-indications and MELD exception awaiting, and their real and predicted survival times in the system.
\\
\\
\\

\paragraph*{ Queue actualization}

The queue actualization at each step of time is handled by the function \textbf{actualization} that acts as follows.
\begin{enumerate}
    \item[] \Step \ $1$: Pick the first non-actualized patient/organ in the queue and get his informations $[a,b_0,b_1,c,d,e]$. If there is no patient/organ left to actualize in the queue go to \Step \ $8$.
    \item[] \Step \ $2$ :  Increase or decrease the relevant information by one time step. (Namely, $b_0,b_1,d$ and $e$ decrease by one  
    and $c$ increases by one). 
    \item[] \Step \ $3$ : 
    If $b_0$ (the real survival time) is negative, the patient/organ is deceased, so we remove the patient/organ from the system and go back to step $1$. 
    If not go to \Step \ $4$.
    \item[]
    \Step \ $4$: If $b_1$ (the predicted survival time) is negative, our prediction was wrong, so we re-draw a new predicted survival time according to $\mu^{'}_{a,c}$.
    \item[]
    \Step \ $5$: If $-1<d<0$, the patient must receive a MELD exception, and we do so using the \textbf{MXP} function, which also updates its information. 
    \item[]
     \Step \ $6$: Using the \textbf{reMELD} function, we determine whether the patient changes MELD class and if he does, we update his 
     information accordingly. We then go back to \Step \ $1$.
    \item[]
    \Step \ $7$: We return the state of the queue, with all patients/organs updated.
\end{enumerate}

\paragraph*{Transition between MELD classes.}

\begin{figure}[h!]
\begin{center}
\includegraphics[scale=0.55]{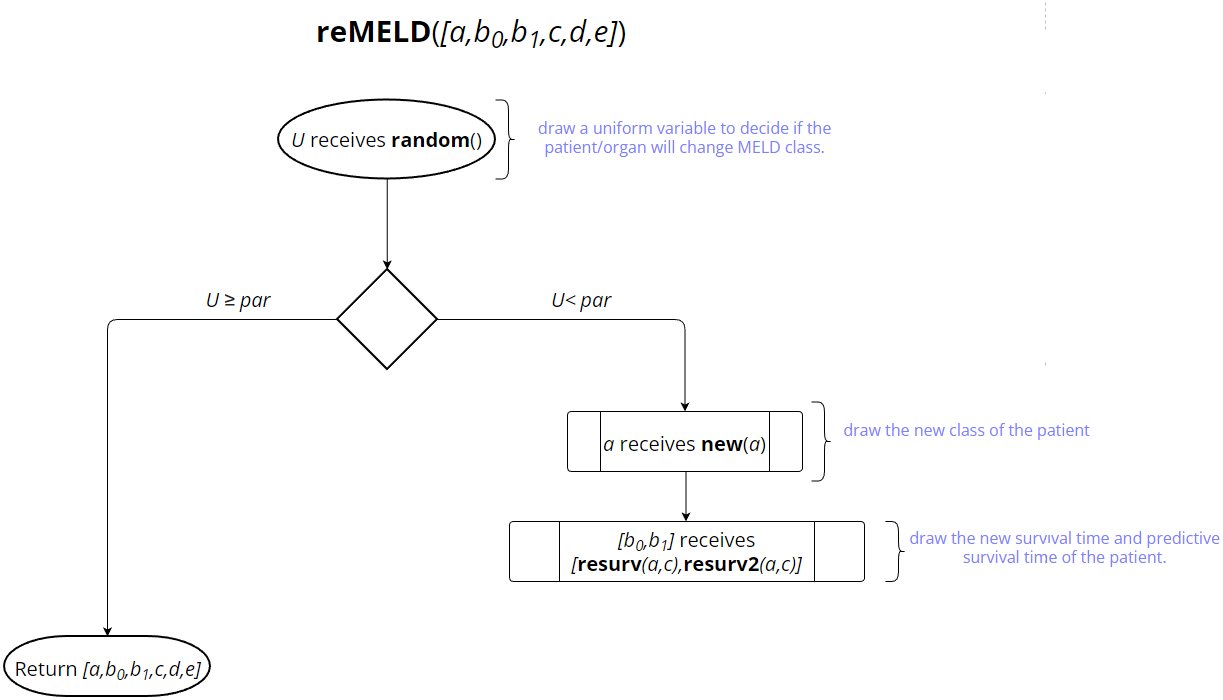}
\caption[smallcaption]{\textbf{reMELD} function}
\label{reMELD}
\end{center}
\end{figure}

The algorithm described in Figure \ref{reMELD} determines randomly if a patient changes MELD class, and if it does, it draws the new class 
according to the $\textit{new}$ function. Last, this routine re-draws the predicted and residual survival times. 


There, the function \textbf{resurv}$(a, c)$ (resp. $\mathbf{resurv_2}(a, c)$)  redraws a remaining real (resp. predicted) survival time according to $\mu_{a,c}$ (resp. $\mu^{'}_{a,c}$), 
and $new(a)$ draws randomly the new MELD of a patient in the following way:
\begin{enumerate}
    \item With probability $\frac{1}{3}$, we draw an improvement of the patient state (lower MELD), and with probability $\frac{2}{3}$, a deterioration of the patient state (higher MELD). (If the MELD is already minimal or maximal, we don't get do anything.) 
\item Then we draw among the possible MELD classes, weighing each class according to the initial arrival probability of each MELD class.
\end{enumerate}

\begin{figure}[h!]
\begin{center}
\includegraphics[scale=0.55]{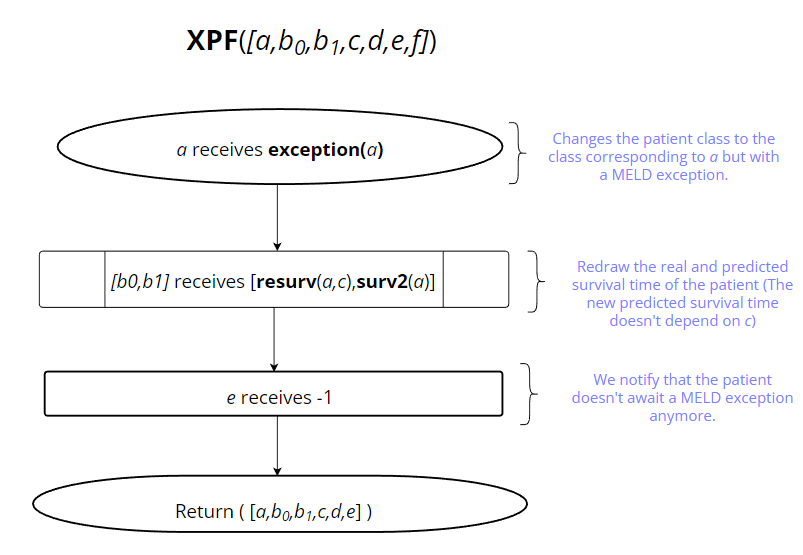}
\caption[smallcaption]{\textbf{XPF} function }
\label{XPF}
\end{center}
\end{figure}

\paragraph*{Assigning a MELD exception}
The program shown in Figure \ref{XPF} updates the class and the real and predicted survival times of a patient that receives a MELD exception.
\\
\\
\\

\bibliographystyle{acm} 
\bibliography{bibliomedicalmath}
\end{document}